\documentclass[12pt, 14paper,reqno]{amsart}
\usepackage{amsmath,amsfonts,amssymb}
\usepackage[breaklinks]{hyperref}
\usepackage{graphicx}

\theoremstyle{plain}
\newtheorem{thm}{Theorem}[section]

\newtheorem{cor}{Corollary}[section]

\newtheorem{thma}{Theorem}

\theoremstyle{proof}

\numberwithin{equation}{section}

\begin{document} 
\title[Class number one problem]{Class number one problem for the real quadratic fields $\mathbb{Q}({\sqrt{m^2+2r}})$}
\author{Azizul Hoque and Srinivas Kotyada}
\address{Azizul Hoque @nstitute of Mathematical Sciences, HBNI, CIT campus, Taramani, Chennai-600113, INDIA}
\email{ahoque.ms@gmail.com}
\address{Srinivas Kotyada @Institute of Mathematical Sciences, HBNI, CIT campus, Taramani, Chennai-600113, INDIA}
\email{srini@imsc.res.in}

\maketitle

\begin{abstract}
We investigate the class number one problem for a parametric family of real quadratic fields of the form $\mathbb{Q}( \sqrt{m^2+4r})$ for certain positive integers $m$ and $r$. 
\end{abstract}

\section{Introduction}
A well-known conjecture of Gauss states that there are infinitely many real quadratic fields with class number one. This is still unresolved. 
Attempts to prove the conjecture have led to important ideas that were instrumental in throwing light in some particular type of quadratic fields, for example for the so-called extended Richaud-Degert type real quadratic fields. Recall that a real quadratic field $\mathbb{Q}(\sqrt{d})$ is of extended Richaud-Degert type if $d$ is square-free positive integer of the form $m^2+r$ with $r\mid 4m$ and $-m<r\leq m$ or $r=\pm 4m/3$. The work of Louboutin \cite{LO1990}, Mollin and Williams \cite{MW1990} and Yokoi \cite{YO1990} confirmed that there are $43$ real quadratic fields, $\mathbb{Q}(\sqrt{d})$; $d = 2, 3, 5, 6, 7, 11, 13, 14, 17, 21, 23, 29, 33, 37, $ $38, 47, 53, 62, 69, 77, 83, 93, 101, 141, 167, 173, 197,  213, 227, 237, 293, 398, 413, 437, 453, 573, ~$   $   677, 717, 1077, 1133, 1253, 1293, 1757$, of extended Richaud-Degert type with class number one with possibly one more such field. However, under the generalized Riemann Hypothesis, there exists at most one more such field with class number one. The problem of finding this exceptional real quadratic field without generalized Riemann Hypothesis is still open.   

Yokoi conjectured in  \cite{YO1986} that there are exactly six real quadratic fields of the form $\mathbb{Q}(\sqrt{m^2+4})$ with class number one, which correspond to $m=1,3,5,7,13,17$. In \cite{BI2003}, Bir\'{o} confirmed this conjecture. Recently, Bir\'{o} and Lapkova \cite{BL2016} extended the result of \cite{BI2003} to a large subclass of Richaud-Degert type real quadratic fields. They proved:
\begin{thma}\label{thmbl}
For odd positive integers $a$ and $m$, let $d = a^2m^2 + 4a$. If $d$ is square-free and $d > 1253$, then $h(d) > 1$.
\end{thma}
Note that by $h(d)$ we mean the class number of the quadratic field $\mathbb{Q}(\sqrt{d})$. As a consequence of Theorem \ref{thmbl}, one can derive the following:
\begin{cor}\label{corbl}
Let $d$ be as in Theorem \ref{thmbl}. Then $5,13,21,29,53,173, 237,293, 437, 453, 1133$ and $1253$ are the only values of $d$ such that $h(d)=1$. 
\end{cor}

The usual method for proving that class numbers of real quadratic fields $K$ in some parametrized families are greater than $1$ is based on continued fractions (cf. \cite{LO1990,LO2000}). This method works only for families for which the continued fraction expansion of the canonical generator $\omega=(d_K+\sqrt{d_K} )/2$ of the ring of algebraic integers of $K$ of discriminant $d_K$ is known beforehand, say for the $\mathbb{Q}(\sqrt{m^2 + 1})$'s with $m^2 + 1$ square-free. 
Another method is based on the computation (in two ways) of special values of zeta function attached to real quadratic fields.  However, this method works only for those fields whose fundamental unit is explicitly known with some other restrictions (cf. \cite{BK1996, CHM2020}). In the present paper we deal with a parametrized family of real quadratic fields for which these methods don't apply. Instead, we will use a modification of the method introduced in \cite[Lemma in p. 218]{ACH65} to prove that some Diophantine equation have no integer solutions. We will prove the following result (which could also be applied to slightly different parametrized families of real quadratic fields):

\begin{thm}\label{thm1}
If $p > 2$ is prime and $d = a^2m^2 +4ap$ with $a > 1$ and $m \geq 1$ is square-free, then the Diophantine equation $x^2- dy^2 = \pm 4p$ has no solution in integers. Hence, $h(d) > 1$.
\end{thm}

\section{Proof of Theorems \ref{thm1}}
Suppose that the equation in Theorem \ref{thm1} has integer solution(s), then without loss of generality, we may assume that $(x,y)$ is a solution with the least possible $y \geq 1$. Then 
\begin{equation}\label{eq2.1}
x^2-dy^2=\pm 4p.
\end{equation} 
Put $\alpha:=\dfrac{x-y\sqrt{d}}{2}$. Then by \eqref{eq2.1}, $N(\alpha)=\pm p$. We now define the following algebraic integer $$\beta:=\dfrac{am^2+2p+m\sqrt{d}}{2}$$ in $\mathbb{Q}(\sqrt{d})$. Then $N(\beta)=p^2$, and 
$$\alpha \beta=\dfrac{(am^2+2p)x-mdy+\left(mx-(am^2+2p)\right)\sqrt{d}}{4}.$$
We now take the norm on both sides and then simplify to get 
\begin{equation}\label{eq2.2}
\pm 4p=\left(\dfrac{(am^2+2p)x-mdy}{2p}\right)^2-\left(\dfrac{mx-(am^2+2p)y}{2p}\right)^2d.
\end{equation}
We can check that both $((am^2+2p)x-mdy)/2p $ and $(mx-(am^2+2p)y)/2p$ are rational integers using the fact $x^2-a^2m^2y=x^2-dy^2+4apy^2$. Thus, 
using the minimality of $y$, we get (from \eqref{eq2.2})
$$y\leq  \left| \dfrac{mx-(am^2+2p)y}{2p}\right|.$$
We first consider the case where
$y\leq (mx-(am^2+2p)y)/2p.$ In this case, we have $(am^2+4p)y\leq mx$, and using this in \eqref{eq2.1} one gets
$\pm 4m^2p\geq\left( (am^2+4p)^2-m^2d\right) y^2$. This implies $\pm m^2p\geq (am^2+4p)py^2$ which further implies $m^2\geq (am^2+4p)y^2$. This is not possible.  

We now consider the remaining case $y\leq ((am^2+2p)y-mx)/2p$. This gives $x\leq amy$. As before, using this in \eqref{eq2.1}, we get
$\pm 4p\leq (am^2-d)y^2$. This implies $\pm p\leq -apy^2$ which gives $1\geq ay^2$. This further implies $a=1$ and $ y=1$, which contradicts the hypothesis.

 We can assume that $\gcd(m, p)=1$; otherwise, it follows from Gauss genus theory that $h(d)>1$. Then $d\equiv a^2m^2\pmod p$, and thus $p$ splits completely in $\mathbb{Q}(\sqrt{d})$ as $(p)=\mathfrak{a}\bar{\mathfrak{a}}$ for some prime ideal $\mathfrak{a}$ and its conjugate $\bar{\mathfrak{a}}$ with norm $N(\mathfrak{a})=p$.
There if $h(d)=1$, then $\mathfrak{a}$ is principal, and hence we can write $\mathfrak{a}=\left(\frac{u+v\sqrt{d}}{2}\right)$ for some $u, v\in \mathbb{Z}$ with $u\equiv v\pmod 2$. This implies that  
$
u^2-dv^2=\pm 4p,
$
which is not possible. This completes the proof. 
\section*{acknowledgements}
The authors are grateful to anonymous referee(s) for careful reading, pointing out a serious error in the previous version and valuable comments which have helped to improve this paper. The authors are also grateful to the referee for drawing the papers \cite{ACH65, LO2000} to their attention.
The authors acknowledge the grants SERB-NPDF (PDF/2017/001958) and SERB MATRICS Project No. MTR/2017/00100,  Govt. of India.

\end{document}